\documentclass{elsarticle}

\usepackage{graphicx} 
\usepackage{amssymb}
\usepackage{xspace}
\usepackage[hyphens]{url}
\usepackage{algorithm}
\usepackage{algpseudocode}
\usepackage{siunitx}

\def\Real#1{{\rm real}(#1)}
\def\Imag#1{{\rm imag}(#1)}
\def\Cmplx#1#2{{#1} + {#2}\, i}
\def\Pain{Painlev{\'{e}}\xspace}
\def\MATLAB{\textsc{Matlab}\xspace}
\def\Eq:ref#1{Eq.~(\ref{#1})}
\def\Alg:ref#1{Algorithm~\ref{#1}}

\providecommand{\keywords}[1]
{
  \small	
  \textbf{\textit{Keywords---}} #1
}

\title{Singularities of the First Painlev{\'{e}} Transcendent}
\author{George F. Corliss, Marquette University, Milwaukee, USA\footnote{Email: {\tt George.Corliss@Marquette.edu}.  Mail: George Corliss, 110 E Sutton Pl, Unit B, Waukesha WI 53188-4041 USA.}}
\date{}

\begin{document}

\begin{abstract}
Consider the solution $y(t)$ for the ordinary differential equation $y' = f(t, y)$ with $t$ complex.  Second-order nonlinear differential equations often exhibit patterns in their poles, branch points, and essential singularities, explored by \Pain and colleagues, 1888--1915.  A variant of the ratio test applied to the Taylor series for the solution $y$ estimates the locations and orders of singularities in the First \Pain Transcendent as an example.  Can you suggest applications in which our singularity location analysis can provide useful insights?
\end{abstract}

\maketitle

\keywords{Painlev{\'{e}} Transcendent; differential equations; complex variables; singularities; Taylor series; seeking applications}

\section{Time is complex?}
\label{sec:TimeIsComplex}

We give numerical maps of the locations of poles of the First \Pain Transcendent, the solution to $y'' = 6 y^2 + t$ in the complex $t$-plane, by extending the work of~\cite{Corliss1980Integrating}.  We explore qualitative views of $y$ enabled by such maps and appeal for applications in which our numerical singularity location can provide useful insights.


The classical theory of ordinary differential equations with complex variables is presented in~\cite{Hartman1964ordinary, Hille1976a, Ince1927a}.  For most equations $y' = f(t, y)$, the solution $y(t)$ is defined for all complex $t$ except for a generally infinite set of points, including poles, branch points, and essential singularities. 

In 1880, Picard~\cite{Picard1880Propriete} first raised questions about second order differential equations with fixed critical points.
In 1885, Poincar{\'e}~\cite{Poincare1885Comptes} showed that differential equations of order two and higher can give rise to previously unknown transcendental functions.  By 1888, Paul \Pain~\cite{Painleve1888lignes}, later joined by Gambier, Fuchs, and Boutroux~\cite{Painleve1895leccons, Painleve1900memoire, borisov2014paul, Boutroux1913a}, studied equations of the form $ y'' = f(t, y, y') $, where $f$ is analytic in $t$ and rational in $y$ and $y'$ (see~\cite{Ince1927a}, Chapter 14, and~\cite{Hille1976a}, pp. 439--447). They characterized equations whose branch points and essential singularities are fixed.  They identified 50 types of such equations, of which 44 are integrable in terms of known functions.  The remaining six types of equations define new transcendental functions.  The simplest such equation is 
  \begin{equation}
    y'' = 6 y^2 + t \, .
    \label{eq:PainOne}
  \end{equation}
Its solution $y$ is the First \Pain Transcendent.  \Pain and his colleagues showed that the First \Pain Transcendent has no branch points or essential singularities.  Its only singularities are poles of order two, whose locations move as continuous functions of the initial conditions.  This paper offers numerical explorations of the locations of poles for a typical set of initial conditions.  We chose \Eq:ref{eq:PainOne} as our illustrative example because of the long and distinguished history of analysis of its singularities. Asymptotic properties are known, but there are no pictures of actual singularity locations.

Section~\ref{sec:RadiusOfConvergence} reviews facts from complex variables and Taylor series.
Section~\ref{sec:ComplexPath} discusses the meaning of a piecewise linear path of integration in the complex plane and the {\tt odets}~\cite{NedialkovPryce2025} algorithm augmented to estimate the locations and orders of singularities.
Section~\ref{sec:PainleveSingularities} gives a picture of singularities of the First \Pain Transcendent.
Section~\ref{sec:BoundaryRowsColumns} fills in the full picture of singularity locations shown in Figure~\ref{fig:PainleveLocations}.
Sections~\ref{sec:ReproducableLocations} through~\ref{sec:SolutionAtTwenty} each pose and answer questions about the behavior of the surface traced by the First \Pain Transcendent from singularity location information.

\medskip
We welcome suggestions of applications for which a study of singularities in the solution to a differential equation can offer insight.

\section{Complex singularities and series radius of convergence}
\label{sec:RadiusOfConvergence}

In this section, we review a few facts from the theory of complex variables and properties of Taylor series that are crucial to the location of singularities.  As an example, consider the function
  \[
    w(z) = \frac{2}{1 + z^2} = \frac{2}{(i - z) (-i - z)} = \frac{i}{i - z}+ \frac{-i}{-i - z}\ , {\rm{\ for\ }} z \ne \pm i \,.
  \]
For all complex $z \ne \pm i$, $w$ is {\it analytic,} which means that all derivatives exist and are finite.
The Taylor series for $w(z)$ expanded at $z_0$,
  \[
    \sum_{k = 0}^\infty w^{(k)}(z_0) (z - z_0)^{k} / k! \, ,
  \]
is represented by a \MATLAB array of Taylor coefficients, 
  \[
    W = \left[ w^{(k-1)}(z_0) (z - z_0)^{k-1} / (k-1)!\right] \, , {\rm {\ for\ }} k = 1, 2, \ldots \, .
  \]
Let $h = z - z_0$.  Then at $z_0 = 0$, 
  \[
    w(z) = 2 \left( 1 + h^2 + h^4 + \ldots\right)
  \]
is a convergent geometric series for $| h | < 1$.  That is, the radius of convergence of the Taylor series for $w(z)$ expanded at $z_0 = 0$ is 1, equal to the distance from the point of expansion $0$ to the closest singularities at $\pm i$.

Generally, for any function $y$, we use $y(t)$ to refer to a single complex value and the capital $Y$ the \MATLAB array of Taylor coefficients at a usually implied point of expansion.

Except for pathological examples, the radius of convergence of a Taylor series of an analytic function is the distance from the point of expansion $z_0$ to the closest singularity (called the primary singularity).  We say that a primary singularity {\it dominates} if no other singularities lie near the circle of convergence.  
The contribution of each secondary singularity drops with the series order as a geometric series in the distance from the point of expansion to each singularity.  Hence, at 30 terms of the series, a secondary singularity whose distance from the point of expansion is twice that of the primary singularity affects only the $30^{\rm th}$ and subsequent binary digits.  
Beyond about 52 terms of the series, double-precision values carry no information whatsoever of secondary singularities more than twice the radius of convergence from the point of expansion.
Computationally, by using a relatively long Taylor series, we see only the effects of a primary singularity unless two or more singularities are nearly the same distance from the point of expansion.  That insight explains why we can estimate locations of primary singularities easily and accurately, as explained in the following sections.

\section{Complex path of integration}
\label{sec:ComplexPath}

Beginning in this section, we switch from $z$ as the independent variable, as is common in complex analysis, to $t$, the customary independent variable in time-dependent differential equations.  To locate singularities, we follow a piecewise linear path of integration in the complex $t$-plane.  For example, in the upper panel of Figure~\ref{fig:PainleveRegion1A}, we locate singularities of $y$ on the real $t$-axis by integrating from $ t = \Cmplx{0}{0}$ to $\Cmplx{1.0}{0.3}$ and then to $\Cmplx{20.5}{0.3}$.

We integrate with the Taylor series ODE package {\tt odets}~\cite{NedialkovPryce2025} with slight modifications, as shown in \Alg:ref{alg:LocSing}.  {\tt Odets} is the same as \Alg:ref{alg:LocSing} except that
  \begin{itemize}
      \item The ``path of integration'' is usually real, for example, $t = 0$ to 20;
      \item Step 2 adds estimates of the location and order of primary singularities.
    \end{itemize}

\begin{algorithm}
\caption{Locate singularities near a path of integration}
\label{alg:LocSing}
\begin{algorithmic}
\Require Piecewise linear integration path in complex t-plane \\
\Loop \ for each integration step $t_k$  \\
1:     \qquad Expand the Taylor series for $y(t)$ at $t_k$;  \\
2:     \qquad Estimate the location and order of the primary singularity
      \\ \qquad \qquad with a measure of confidence;  \\
3:     \qquad Estimate a step to $t_{k+1}$ to satisfy absolute and relative error criteria;  \\
4:     \qquad Sum the series to get $y(t_{k+1})$ and $y'(t_{k+1})$, initial conditions 
      \\ \qquad \qquad for the next step;
\EndLoop \\ \\
\Return \ Steps, solution values, estimated singularity locations, orders, and confidences
\end{algorithmic}
\end{algorithm}

To specify a complex line segment to {\tt odets}, we give a complex {\tt tspan = [t0, tf]}, and {\tt odets} integrates along the straight line segment from {\tt t0} to {\tt tf}.  For example, the script
{\small
\begin{verbatim}
    pain1ode = @(t, y) [ y(2); 6 * y(1)^2 + t ];
    t0 = 0 + 0i;
    y0 = [ 0.5 + 0.0i, 0.9 + 0.0i ];
    tspan = [ t0, 1.0 + 0.3i ];
    opts = odetsset(AbsTol = 1e-10, RelTol = 1e-10, TSOrder = 45);
    [t, y] = odets(pain1ode, tspan, y0, opts);
\end{verbatim}}
integrates from $ t = \Cmplx{0}{0}$ to $\Cmplx{1.0}{0.3}$ in Figure~\ref{fig:PainleveRegion1A}.  With suitable formatting, this script produces
{\small
\begin{verbatim}
 Step   Complex time t                      y                     y'
    1  0.000 + 0.000 i    0.5000 +   0.0000 i    0.9000 +   0.0000 i
    2  0.500 + 0.150 i    1.2357 +   0.4654 i    2.7278 +   1.5701 i
    3  0.750 + 0.225 i    1.8174 +   1.6192 i    3.6538 +   6.7531 i
    4  1.000 + 0.300 i    0.6948 +   4.6897 i  -10.8694 +  17.4080 i
\end{verbatim}}
showing $y$ and $y'$ at the endpoints and at two intermediate steps.  Commonly used \MATLAB solvers such as {\tt ode45} have similar capabilities to follow a linear path of integration in the complex plane.  We use {\tt odets} because it generates a Taylor series that we use for singularity location.

In \MATLAB, a function that works with doubles generally works with complex numbers.  In the version of {\tt odets} we received, the only edit required was to change transpose ({\tt '}) to element-wise transpose ({\tt .'}) in two places.  Then to {\tt odets}, we added a function call to estimate the location and order of the primary singularity with a measure of confidence.

\section{Singularities of the First \Pain Transcendent}
\label{sec:PainleveSingularities}

Having software {\tt odets} to follow a piecewise linear path of integration in the complex $t$-plane, consider the initial value problem 
  \begin{equation}
    y'' = 6 y^2 + t;\ y(0) = 0.5,\ y'(0) = 0.9\, .
    \label{Eq:IVP}
  \end{equation}
For our purposes, the initial conditions are not distinguished; other initial conditions yield similar pictures.

\begin{figure}[ht]
  \vspace{-10 pt}
  \begin{center}
  \includegraphics[width=\linewidth]{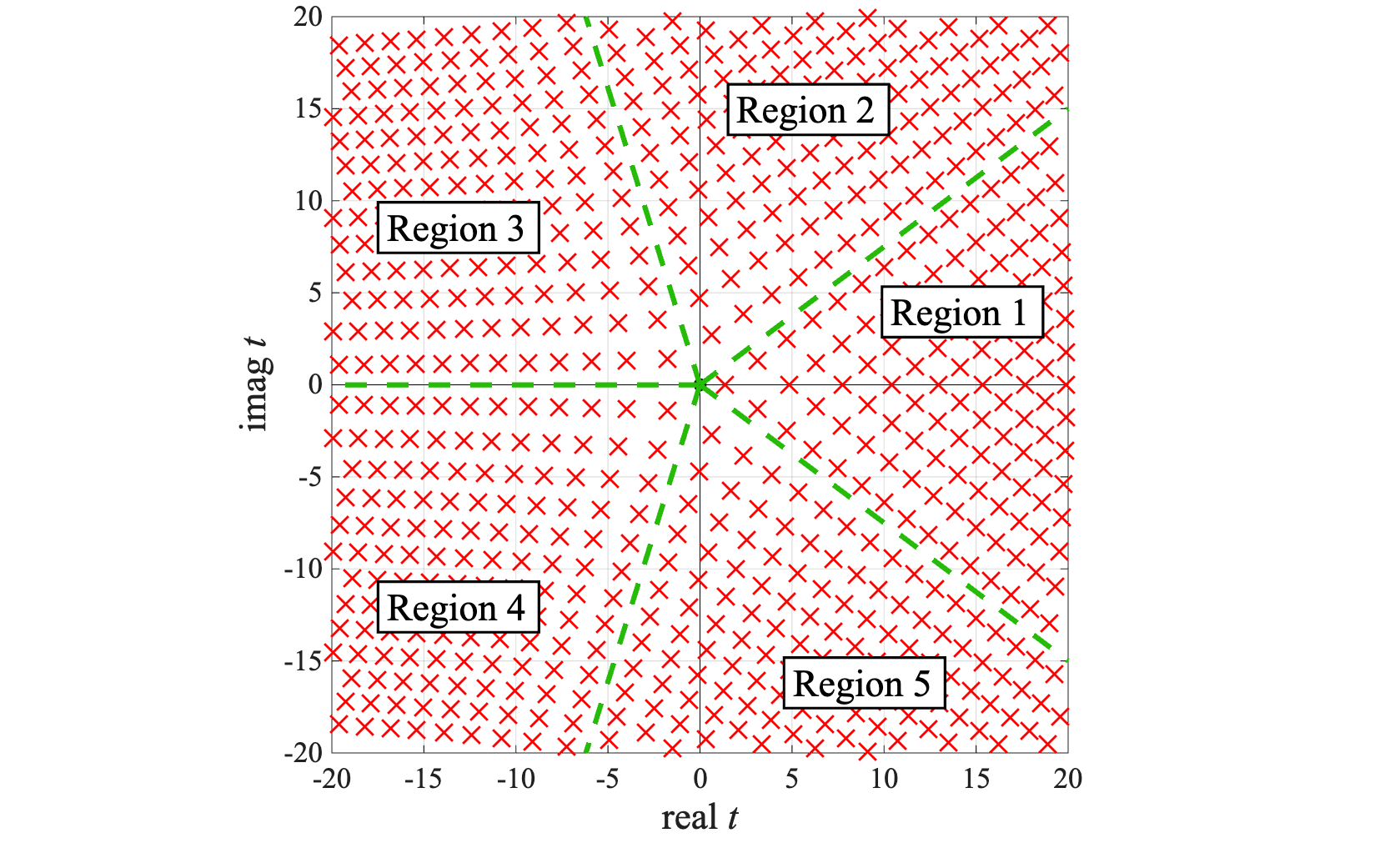}
  \end{center}
  \vspace{-15 pt}
  \caption{Singularity locations for First \Pain Transcendent with $y(0) = 0.5$, $y'(0) = 0.9$.}
  \label{fig:PainleveLocations}
\end{figure}

Since the solution $y$, the First \Pain Transcendent, is real-valued on the real $t$-axis, any singularities are real, or they appear in conjugate pairs.  Figure~\ref{fig:PainleveLocations} shows the singularities of $y$ for $-20 \le \Real{t} \le 20$, $-20 \le \Imag{t} \le 20$.  In the rest of this paper, we study only the upper half of the complex $t$-plane.  This section and the next describe how Figure~\ref{fig:PainleveLocations} was generated.  

We estimate the locations and orders of singularities using the Three-Term Analysis of Chang and Corliss~\cite{ChangCorliss1980ratio}, an extension of the ratio test for convergence of a Taylor series.  In the neighborhood of a primary singularity at $ t = a $ of order $s$, the terms of the Taylor series for $y$ the solution expanded at $ t = 0 $ are asymptotic to the terms of the Taylor series for the model problem
  \[
      w(t) = (a -t)^{-s} \, .
      \label{eq:Wmodel}
  \]
For the First \Pain Transcendent, $s = 2$ was established by \Pain, but we retain the general case since the value of $s$ is generally not known in advance. If $x_1 = 1/a$ and $x_2 = s/a$, the terms of the Taylor series for $w(t)$ satisfy
  \[
      k \, W_{k+1} = (k-1) \, W_k \, x_1 + W_k \, x_2 \, , {\rm \ with\ } W_1 = a^{-s}, \ W_2 = s \, a^{-s-1} \, .
      \label{eq:Wrecur}
  \]
Suppose that we have computed $p$ terms of $Y$, the Taylor series for $y$ expanded at the current integration step.  We use $p = 45$.  Long Taylor series allow the ODE integration algorithm to take integration steps much longer than usual methods such as {\tt ode45} and produce singularity locations that are more accurate than those from shorter series because secondary singularities have less influence on terms in the tail of the longer series. Automatic differentiation allows {\tt odets} to compute long Taylor series accurately and efficiently.  We solve the linear system 
  \begin{equation}
      k \, Y_{k+1} = (k-1) \, Y_k \, x_1 + Y_k \, x_2 \, , {\rm{\ for\ }} k = p-6,\, \ldots,\, p-1 \, ,
      \label{eq:Yrecur}
  \end{equation}
for unknowns $x_1$ and $x_2$, let $a = 1 / x_1$ and $s = x_2 / x_1$, and compute a relative residual.  Using $p-6$ is chosen heuristically.  If we use two copies of \Eq:ref{eq:Yrecur} with $k = p-2$ and $p-1$, we find a solution, and the residual is zero, whether the series for $w(t)$ is a good model for the tail of the series for $y(t)$ or not.  By using $k$ from $p-6$ to $p-1$, we have six equations in two unknowns.  If the relative residual is quite small (we use a heuristically-determined acceptance threshold of $10^{-10}$), the tail of the series for $y(t)$ is well-approximated by the tail of the series for $w(t)$, and the estimates for the singularity location at $t_k + a$ and the singularity order $s$ are good, typically better than single precision, and we accept them.  If the linear system \Eq:ref{eq:Yrecur} is ill-conditioned or the relative residual is not small, $w$ is not a good model for $y$, and we reject the singularity location and order estimates.

\section{Locating boundary rows and columns of singularities}
\label{sec:BoundaryRowsColumns}

In Figure~\ref{fig:PainleveLocations}, we observe three distinct wedge-shaped regions, approximately centered on the fifth roots of unity ($\Cmplx{\cos(2 \pi k / 5)}{\sin(2 \pi k / 5)}$, for $ k = 0, 1, \ldots, 4$), as suggested by \Pain.  We locate singularities in each region separately.

\begin{figure}[ht]
  \begin{center}
  \includegraphics[width=\linewidth]{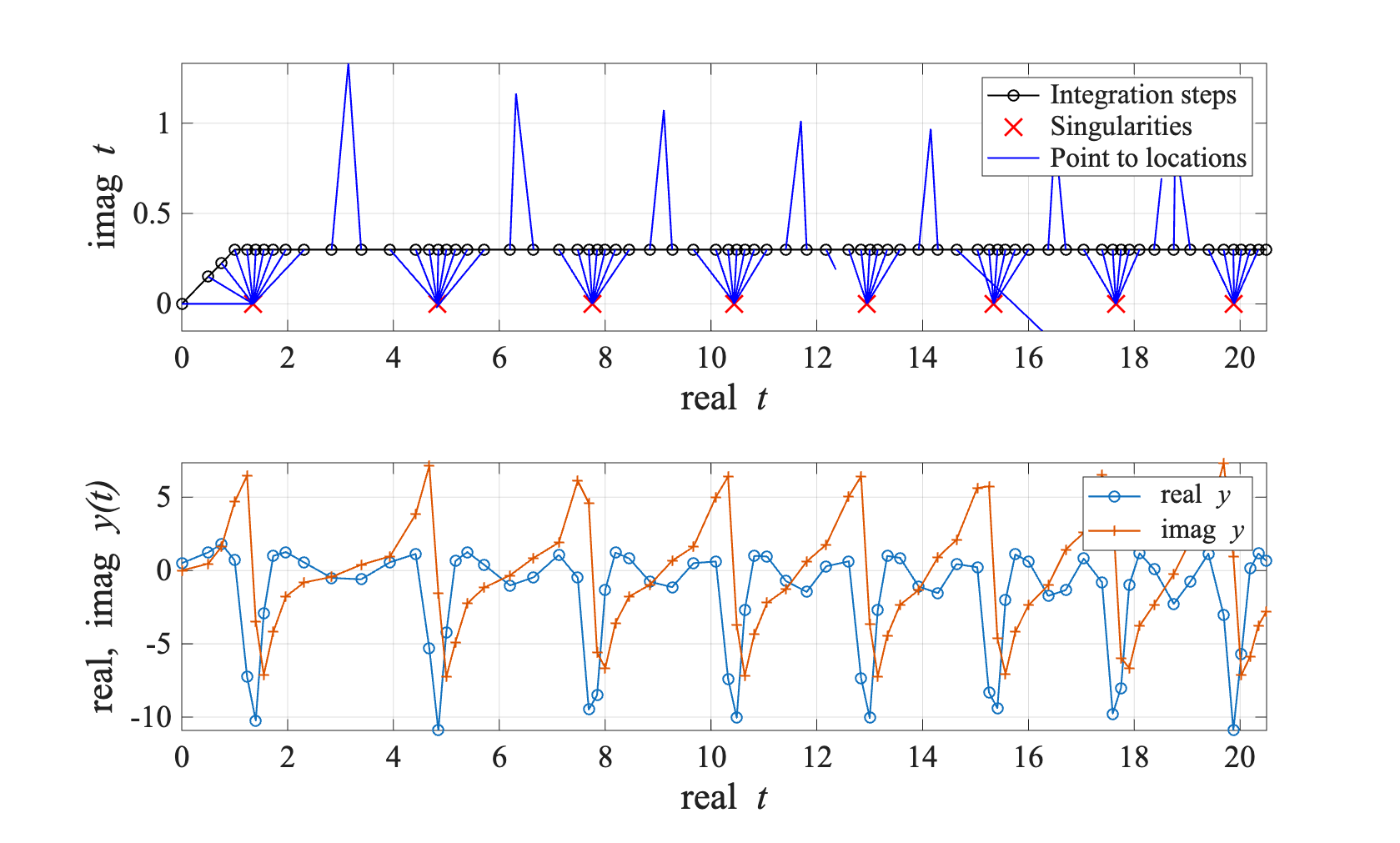}
  \end{center}
  \vspace{-15 pt}
  \caption{Real singularities for the First \Pain Transcendent. Upper panel shows the path of integration and the locations of real-valued singularities.  Lower panel shows real and imaginary components of the solution along the path of integration.}
  \label{fig:PainleveRegion1A}
\end{figure}

First, we locate the singularities on the positive real $t$-axis applying \Alg:ref{alg:LocSing} along the path $\Cmplx{0}{0}$ to $\Cmplx{1}{0.3}$ to $\Cmplx{21}{0.3}$, as shown in Figure~\ref{fig:PainleveRegion1A}.  For convenience, we wrote a wrapper function to accept an array of complex points, call {\tt odets} as above for each line segment, and marshal the solution at each step from the multiple segments.  The script
{\small
\begin{verbatim}
    pain1ode = @(t, y) [ y(2); 6 * y(1)^2 + t ];
    t0 = 0.0 + 0.0i;
    y0 = [ 0.5 + 0.0i, 0.9 + 0.0i ];
    tspan = [ t0, 1.0 + 0.3i, 20.5 + 0.3i ];
    opts = odetsset(AbsTol = 1e-10, RelTol = 1e-10, TSOrder = 45);
    [ t, y, te, ye, ie, location, order, confidence ]  ...
            = odets_multiSegment(pain1ode, tspan, y0, opts);
\end{verbatim}}
integrates from $ t = \Cmplx{0}{0}$ to $ \Cmplx{1.0}{0.3}$ to $ \Cmplx{20.5}{0.3}$ in Figure~\ref{fig:PainleveRegion1A}.  With suitable formatting, this script produces
{\small
\begin{verbatim}
Step      Complex time t                      y                     y'
   1   0.000 +   0.000 i    0.5000 +   0.0000 i    0.9000 +   0.0000 i
   2   0.500 +   0.150 i    1.2357 +   0.4654 i    2.7278 +   1.5701 i
   3   0.750 +   0.225 i    1.8174 +   1.6192 i    3.6538 +   6.7531 i
   4   1.000 +   0.300 i    0.6948 +   4.6897 i  -10.8694 +  17.4080 i
   5   1.234 +   0.300 i   -7.2441 +   6.4638 i  -53.5809 + -27.9573 i
    . . .
  73  20.345 +   0.300 i    1.1646 +  -3.7448 i   -0.2263 +   9.5219 i
  74  20.500 +   0.300 i    0.6431 +  -2.7799 i   -5.6872 +   3.7117 i
\end{verbatim}}
showing $y$ and $y'$ agreeing with the plots for $y(t)$ in the lower panel of Figure~\ref{fig:PainleveRegion1A}.

In the upper panel of Figure~\ref{fig:PainleveRegion1A}, the line of small black circles represents the integration steps of {\tt odets} Step 1 in \Alg:ref{alg:LocSing}.  The red $\times$s are singularities. The density of singularities increases slightly with distance from the origin.  The blue segments join integration steps with the singularity location estimated by \Alg:ref{alg:LocSing} at that step.  When several blue segments appear to end at a red $\times$, the location estimates from several integration steps agree, raising our level of confidence in the estimates.

In  Figure~\ref{fig:PainleveRegion1A}, when the integration steps are between the real singularities, \Alg:ref{alg:LocSing} estimates singularities from the next row up (see Figure~\ref{fig:PainleveLocations}), but only for two or three steps, and the estimates of confidence are low.  Hence, we do not count those as ``located,'' and there are no red $\times$s.  At steps near $t = 12$ and $16$, the locations estimated by \Alg:ref{alg:LocSing} were bad, but the estimated confidence was low, and other estimates did not confirm the location estimate.

\medskip
In general, if the primary singularity is dominant, the estimates are good, they agree on successive integration steps, and the estimate of confidence is high.  If there are secondary singularities slightly outside the circle of convergence, the estimates of $a$ and $s$ are poor and inconsistent, the estimate of confidence is low, we discard the estimates, and we do not mark them with red $\times$s.

In the lower panel of Figure~\ref{fig:PainleveRegion1A}, we plot the real (blue) and the imaginary (orange) components of the solution $y(t)$ vs.\hspace{-2 pt} $\Real{t}$.  Toward the left, the pattern differs because we are on the diagonal path of integration.  Once we begin following the horizontal path of integration, the solution appears nearly periodic, although the ``period'' shortens as the real singularities become a little closer together with increasing distance from the origin.

\medskip
At this point, we should be developing a sense of the behavior of the solution $y$, and we should understand how \Alg:ref{alg:LocSing} estimates the locations of the row of singularities on the real $t$-axis shown in Figure~\ref{fig:PainleveRegion1A}.  We proceed to fill in additional singularity locations shown in Figure~\ref{fig:PainleveLocations}.

\begin{figure}[ht!]
  \vspace{-10 pt}
  \begin{center}
  \includegraphics[width=\linewidth]{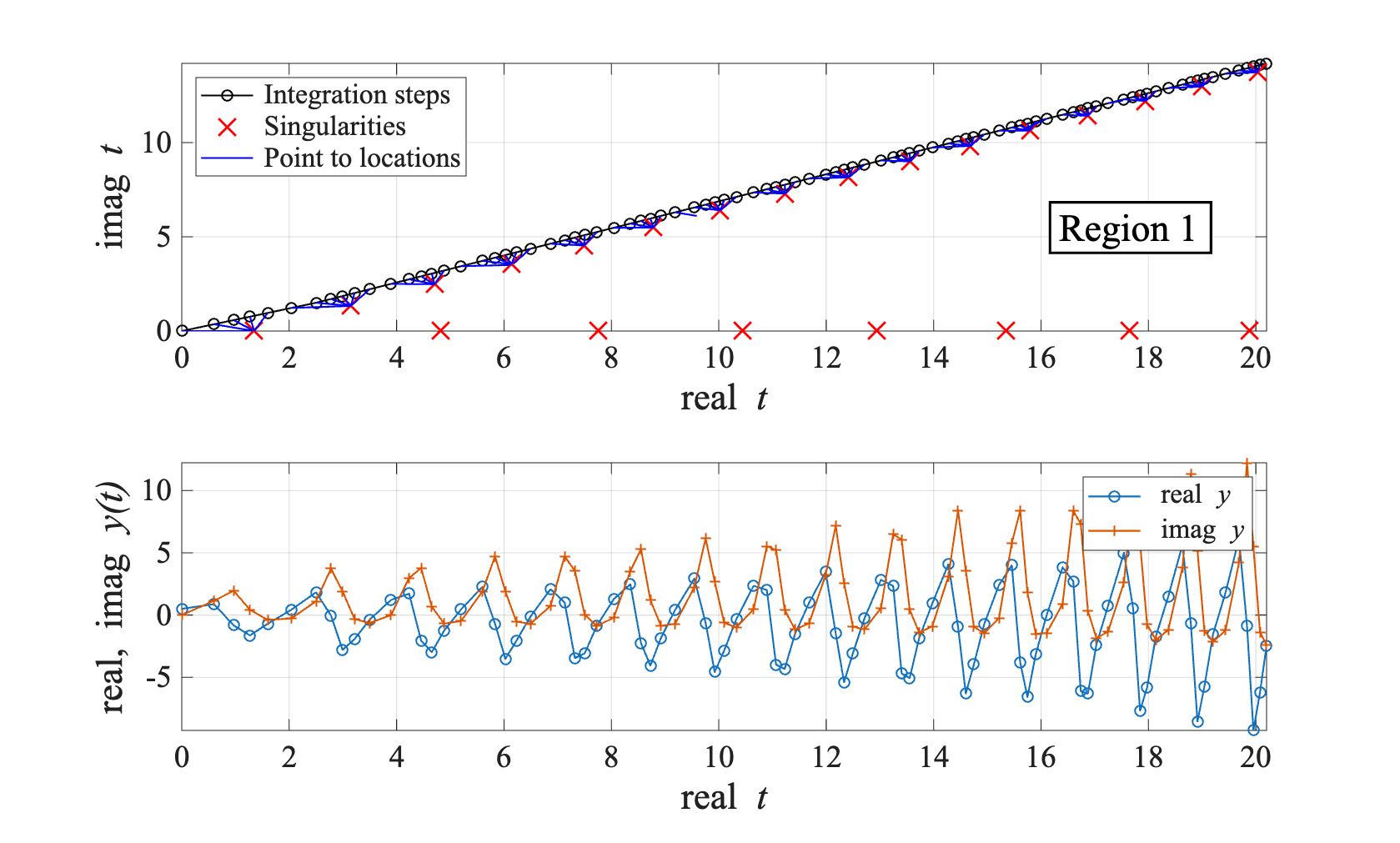}
  \end{center}
  \vspace{-20 pt}
  \caption{Upper panel shows the path of integration and the locations of singularities forming the upper edge of Region 1.  Lower panel shows real and imaginary components of the solution along the path of integration.}
  \label{fig:PainleveRegion1B}
\end{figure}

Next in Figure~\ref{fig:PainleveRegion1B}, we estimate locations of singularities forming a diagonal heading northeastward from the origin in Figure~\ref{fig:PainleveLocations}.  We integrate along a path closely paralleling the desired singularities.  The path must be close enough to the desired singularities to avoid ``seeing'' secondary singularities and far enough away from the desired singularities to allow the step size algorithm to make progress.  Finding this path required substantial trial and error.

The lower panel in Figure~\ref{fig:PainleveRegion1B} shows the real and imaginary components of $y(t)$ vs. $\Real{t}$ along the path of integration, indicated by the line of small black circles in the upper panel.  It appears that the solution becomes more oscillatory as the singularities become closer together, although this is also a consequence of the path approaching each singularity more closely.

\begin{figure}[ht]
  \begin{center}
  \includegraphics[width=0.8\linewidth]{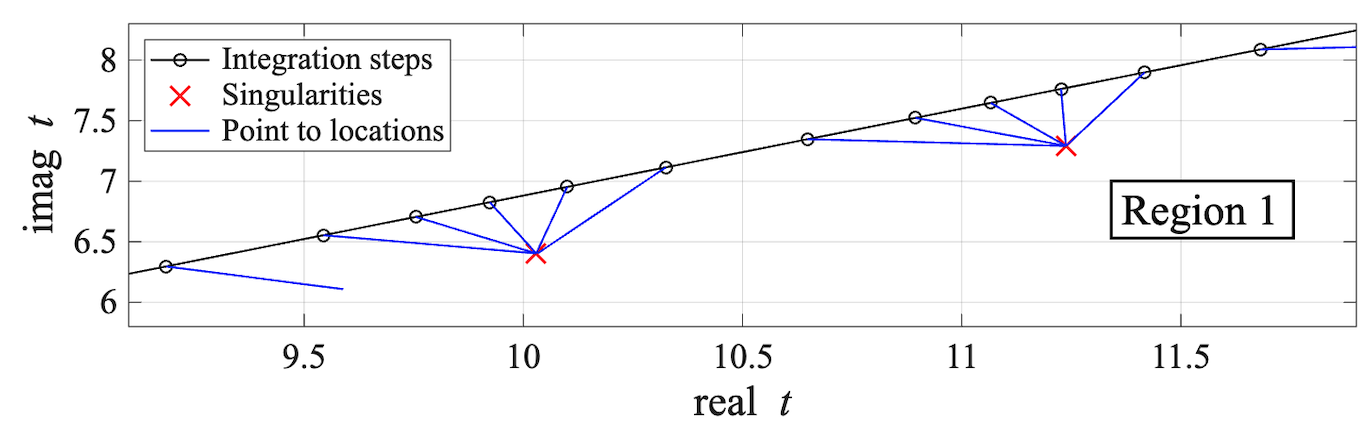}
  \end{center}
  \vspace{-15 pt}
  \caption{Detail of Figure~\ref{fig:PainleveRegion1B} near $t = \Cmplx{10}{7}$ showing details of location of two singularities.  Location estimates show excellent agreement.}
  \label{fig:PainleveRegion1Bz}
\end{figure}

From Figure~\ref{fig:PainleveRegion1B}, we zoom into a neighborhood of $t = \Cmplx{10}{7}$ in Figure~\ref{fig:PainleveRegion1Bz}.  We see that the path of integration is close enough to each singularity for robust performance of the singularity location in \Alg:ref{alg:LocSing}.

\begin{figure}[ht!]
  \vspace{-10 pt}
  \begin{center}
  \includegraphics[width=\linewidth]{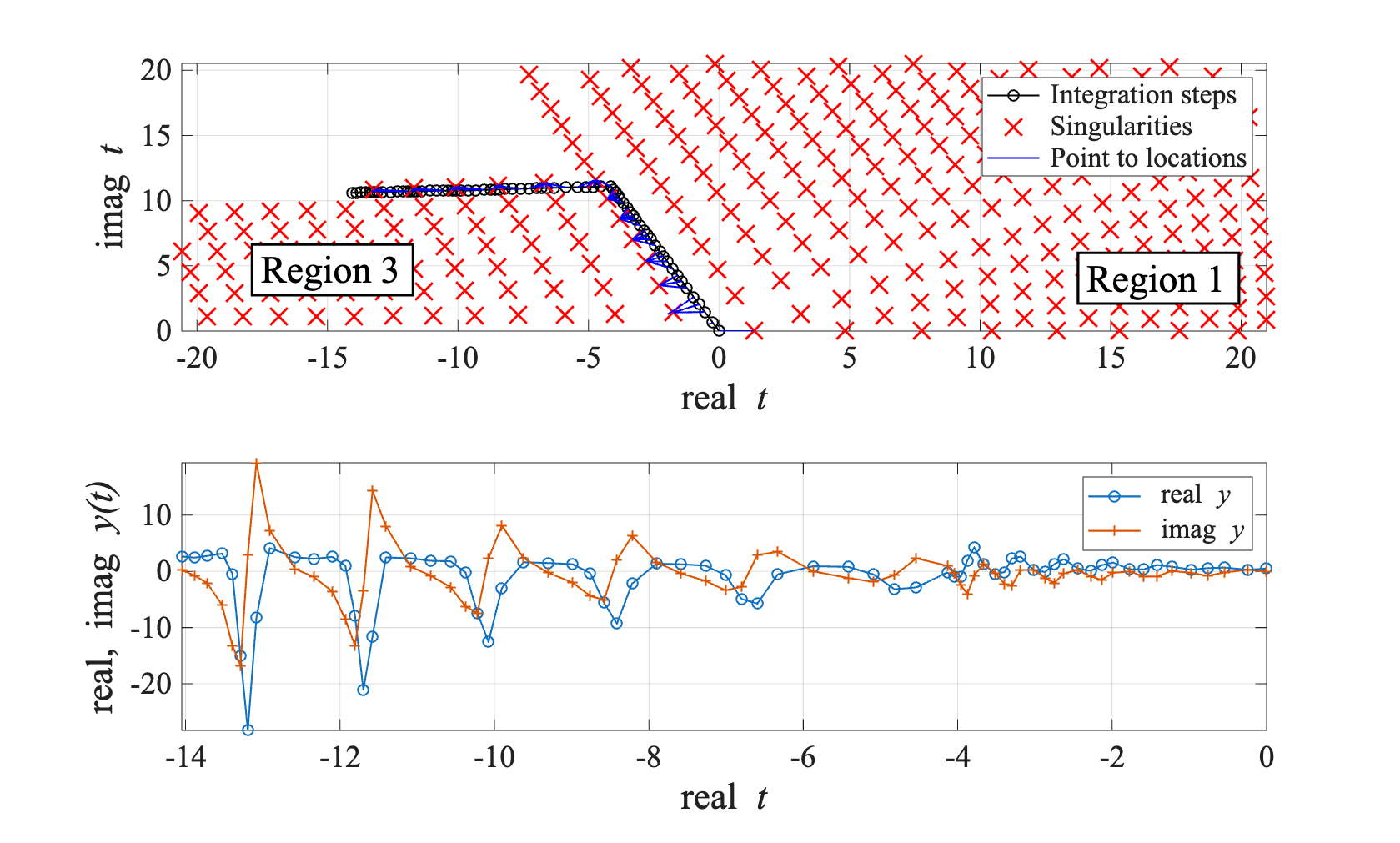}
  \end{center}
  \vspace{-20 pt}
  \caption{Upper panel shows previously located singularities and the path of integration to discover several singularities in the interior of Region 3.  Lower panel shows real and imaginary components of the solution along the path of integration.}
  \label{fig:PainleveFilling395}
\end{figure}

In a similar manner, we locate singularities forming the edges of singularity Regions 2 and 3 seen in Figure~\ref{fig:PainleveLocations}.  Considerable trial and error is required to find suitable paths of integration, as singularities get closer together with increasing distance from the origin.  Also, the singularities do not lie on straight lines.  In Figure~\ref{fig:PainleveFilling395}, Regions 1 and 2 have been filled by integrating along suitable piecewise linear paths, and filling Region 3 is partially completed.  Suitable paths of integration are built dynamically as each new singularity location is acquired until we have the singularity locations shown in Figure~\ref{fig:PainleveLocations}.

\begin{figure}[ht!]
  \begin{center}
  \includegraphics[width=0.8\linewidth]{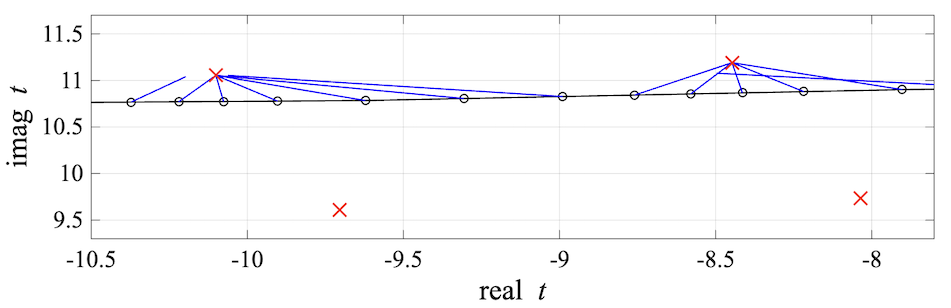}
  \end{center}
  \vspace{-15 pt}
  \caption{Detail of Figure~\ref{fig:PainleveFilling395} near $t = \Cmplx{-9}{10}$ in Region 3 showing details of location of two singularities. Location estimates are in good agreement.}
  \label{fig:PainleveNearm9}
\end{figure}

\newpage

\section{Are locations reproducible?}
\label{sec:ReproducableLocations}

The singularities shown in Figures~\ref{fig:PainleveLocations} --~\ref{fig:PainleveNearm9} are located by analytic continuation of the solution $y$ using 45-term Taylor series along relatively long paths of integration in the complex $t$-plane.  We are taking several tens of integration steps passing close to many singularities, all using standard IEEE double-precision arithmetic in \MATLAB.  We should be surprised if {\bf anything} meaningful survives the accumulation of round-off error, but it does!

\begin{figure}[ht]
  \vspace{-10 pt}
  \begin{center}
  \includegraphics[width=\linewidth]{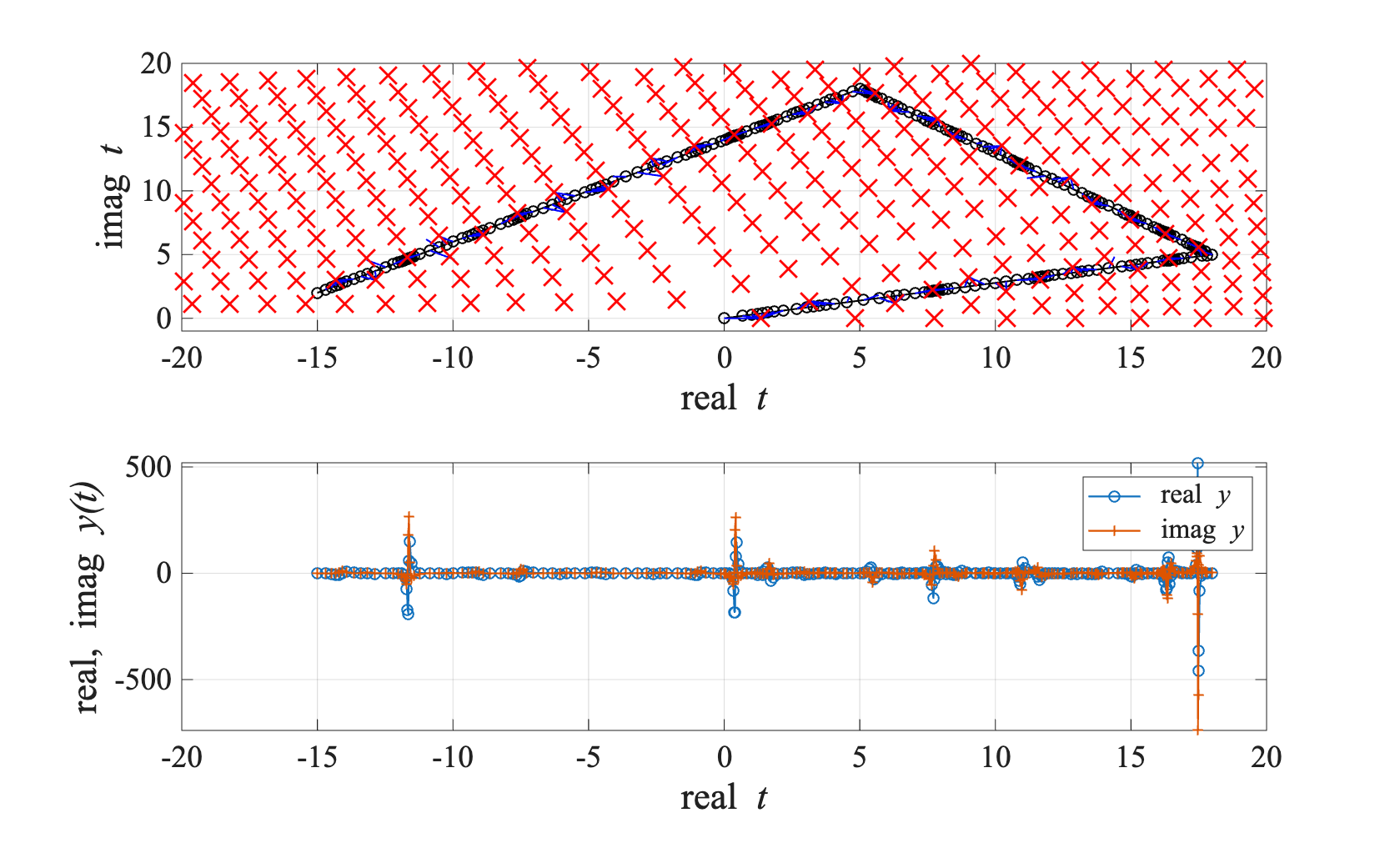}
  \end{center}
  \vspace{-10 pt}
  \caption{Locations do {\bf not} depend on the path of integration!
  Upper panel shows previously identified singularities and a path of integration from the origin left to $\Cmplx{18}{5}$, back in the negative direction to $\Cmplx{5}{18}$, and ending at $\Cmplx{-15}{2}$.  Lower panel shows real and imaginary components of the solution along the path of integration.  Relatively large solution values occur when the path is very close to a singularity.
  }
  \label{fig:PainleveWander}
\end{figure}

\begin{figure}[ht]
  \begin{center}
  \includegraphics[width=0.8\linewidth]{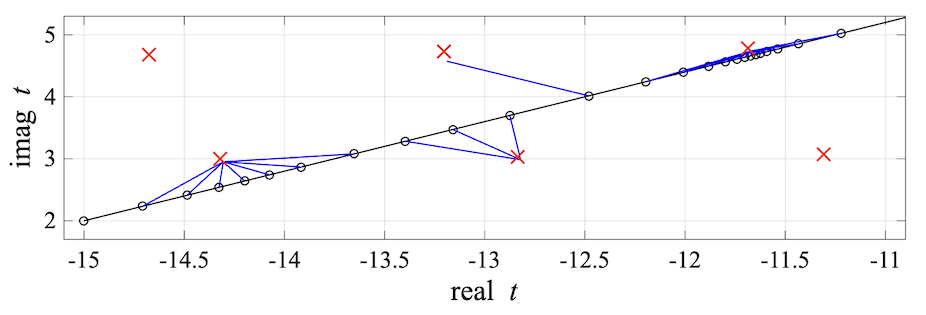}
  \end{center}
  \vspace{-15 pt}
  \caption{Detail of Figure~\ref{fig:PainleveWander}  near $t = \Cmplx{-13}{3}$ showing good agreement with previously identified singularities.}
  \label{fig:PainleveNearm13}
\end{figure}

If we follow a different path of integration, do we encounter the same singularities?  Yes, the singularity locations are stable with respect to the path of integration, as suggested by Figure~\ref{fig:PainleveWander}. The piecewise-linear path has vertices $\Cmplx{0}{0}$, $\Cmplx{18}{5}$, $\Cmplx{5}{18}$, and $\Cmplx{-15}{2}$.  We followed a path for 264 integration steps and located 36 previously identified singularities.  Figure~\ref{fig:PainleveNearm13} shows that even at the end of such a journey, \Alg:ref{alg:LocSing} ``sees'' the same singularities.

The lower panel of Figure~\ref{fig:PainleveWander} merits some explanation.  Like the similar panel in Figure~\ref{fig:PainleveRegion1A}, we plot the real (blue) and imaginary (orange) parts of the solution $y(t)$ vs. $\Real{t}$ for the path of integration shown in the upper panel.  Since the integration path starts at the origin, goes left to $\Cmplx{18}{5}$, turns to head in the negative direction to $\Cmplx{5}{18}$, and ends at $\Cmplx{-15}{2}$, in the lower panel, the solution graphs for $\Real{t} >= 0$ are doubled.  The extremes occur where the path of integration approaches very close to singularities, testing the stability of subsequent locations.

This section has offered evidence that the singularity locations shown in Figure~\ref{fig:PainleveLocations} are independent of the integration paths used to find them.  That is not the case for a similar analysis of the solution to the Van der Pol oscillator equation, whose singularities are branch points, not poles (work in progress).

\section{Closed loop path of integration}
\label{sec:ClosedLoopPath}

As a second check on the stability of singularity locations with respect to the path of integration, we follow a closed loop path of integration.  The computed solution returns to the original initial conditions, as expected.

\begin{figure}[ht]
  \begin{center}
  \includegraphics[width=0.9\linewidth]{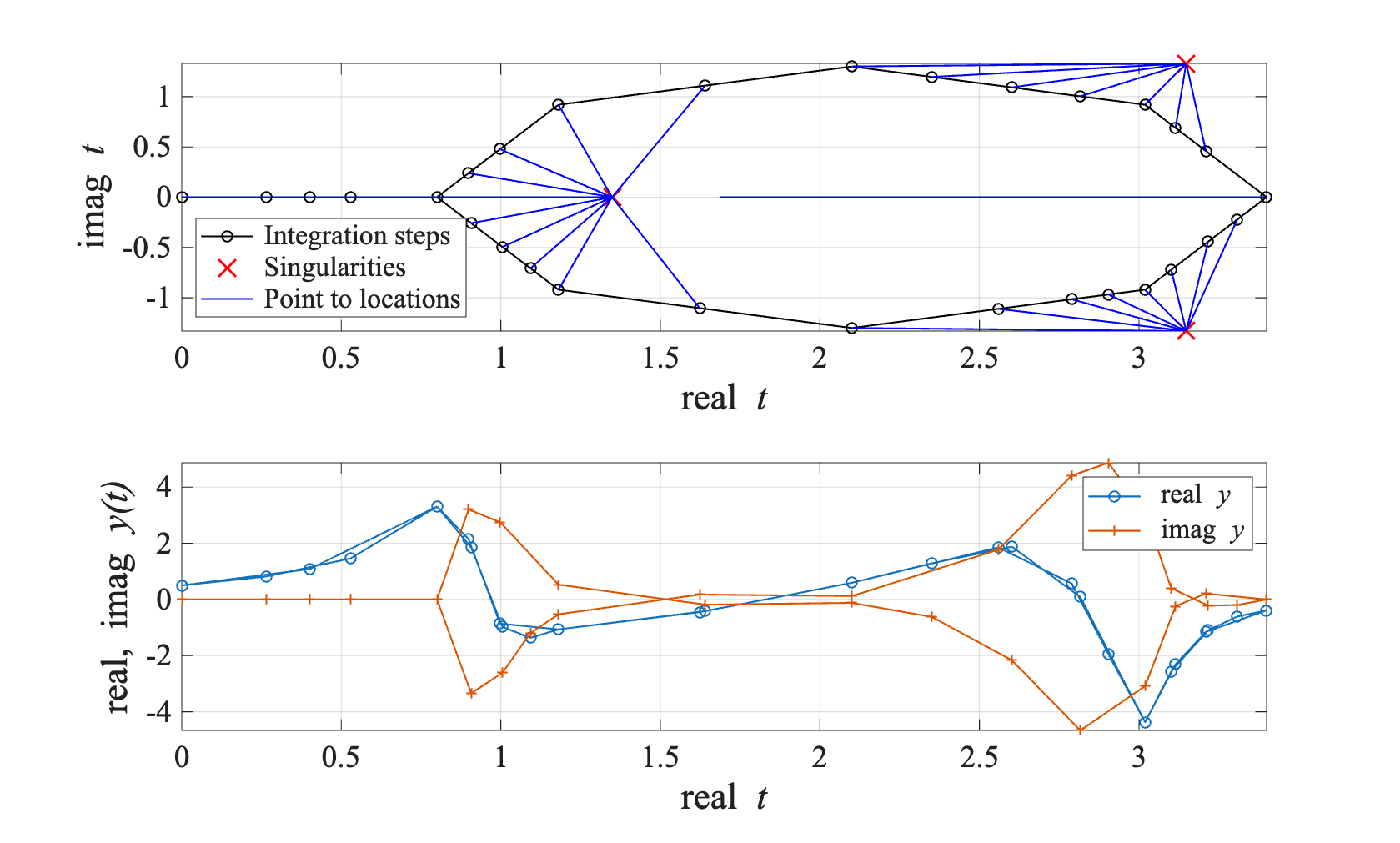}
  \end{center}
  \vspace{-15 pt}
  \caption{Upper panel shows a path of integration from the origin to $\Cmplx{0.8}{0}$, counterclockwise around the hexagon, and returning to the origin.  Integration encounters and located singularities at $\Cmplx{1.3486}{0}$ and $3.1480 \pm 1.3306\, i$.  Lower panel shows that real and imaginary components of the solution return to its starting value, as expected.}
  \label{fig:PainleveCirclePole}
\end{figure}

In Figure~\ref{fig:PainleveCirclePole}, we begin at the origin with the same initial conditions $y(0) = 0.5$ and $y'(0) = 0.9$.  We previously determined that $y(t)$ has a singularity (pole of order 2) at $t \approx \Cmplx{1.3486}{0}$.  The integration path visits $t = \Cmplx{0.8}{0}$, follows an elongated hexagon in the counterclockwise direction around the singularity, and returns to the origin.  The hexagon could be centered on the singularity, but we elongated it so that several integration steps from \Alg:ref{alg:LocSing} locate the conjugate pair of singularities at $t \approx 3.1480 \pm 1.3306\, i$, just to add interest.

On return to the origin, the computed solution is
  \[
    [y(0), y'(0)] = [0.5000000065 + \num{1.77e-9}\, i, 0.8999999889 + \num{3.51e-9}\, i]\ ,
  \]
compared with the true $[0.5, 0.9]$.  The Euclidean error is $[0.67, 1.17] \times 10^{-8}$, while the solution was computed using {\tt AbsTol = 1e-10, RelTol = 1e-10}.

The lower panel in Figure~\ref{fig:PainleveCirclePole} merits some explanation.  Like the similar panel in Figure~\ref{fig:PainleveRegion1A}, we plot the real (blue) and the imaginary (orange) components of the solution $y(t)$ vs. $\Real{t}$ for the path of integration shown in the upper panel.  Since the path of integration starts at the origin, goes counterclockwise around the pole, and returns to the origin, the real and the imaginary components of $y(t)$ track right and then left.  The real component (blue) tracks the same in both directions, while for the imaginary component (orange), the two directions are complex conjugates.

What is important in the lower panel of Figure~\ref{fig:PainleveCirclePole} is that at $t = 0$, the starting and ending points are indistinguishable, differing by less than $\num{7e-9}$. 

Yes, we return to the original initial conditions.  After 31 integration steps with absolute and relative accuracy tolerances $\num{1e-10}$ and encountering three singularities, the computed values of $y(0)$ agree with the initial conditions with an error less than $\num{7e-9}$.  This check raises our level of confidence in the accuracy of integration around a pole.

\section{Qualitative solution surface}
\label{sec:QualitativeSolution}

This section was inspired by a paragraph in my old calculus book by Fobes and Smyth~\cite{Fobes1963Calculus} in which they described a motley group of tourists walking around the surface of a hyperbolic paraboloid near its saddle point.  The authors invited us to visualize a surface as a landscape to expose visualizations of solutions to differential equations enabled by our \Alg:ref{alg:LocSing}.

Both the independent variable $t$ and the dependent variable $y$ are complex, so the solution ``surface'' is 4-dimensional.  Here, we choose to think of the components $\Real{y}$ and $\Imag{y}$ vs.\hspace{-3 pt} complex $t$ as if they are separate, related 3-dimensional surfaces, as if we are walking in the 2-dimensional $t$-plane independently exploring the $\Real{y}$ surface and the $\Imag{y}$ surface.

\begin{figure}[ht!]
  \vspace{-10 pt}
  \begin{center}
  \includegraphics[width=0.98\linewidth]{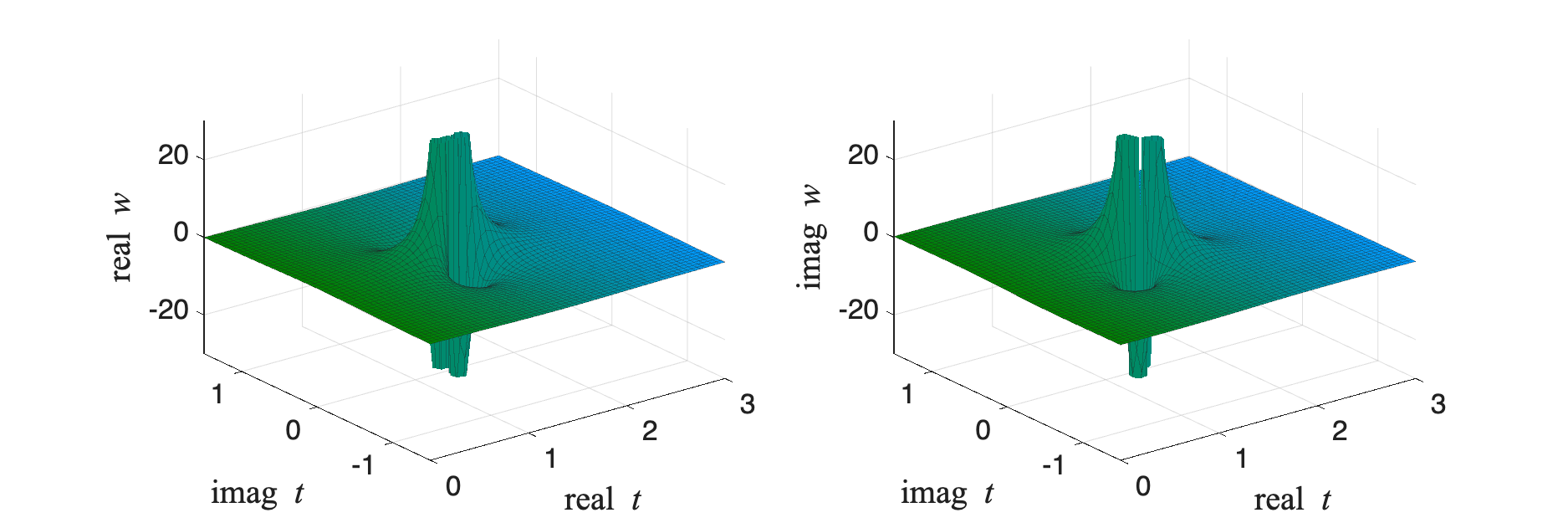}
  \end{center}
  \vspace{-15 pt}
  \caption{Real and imaginary parts of $w(t) = (a - t)^{-2}$.  Each has two spikes and two pits. The First \Pain Transcendent is similar in the neighborhood of each of its poles.}
  \label{fig:PainleveRealImag2}
\end{figure}

\begin{figure}[ht]
  \begin{center}
  \includegraphics[width=0.9\linewidth]{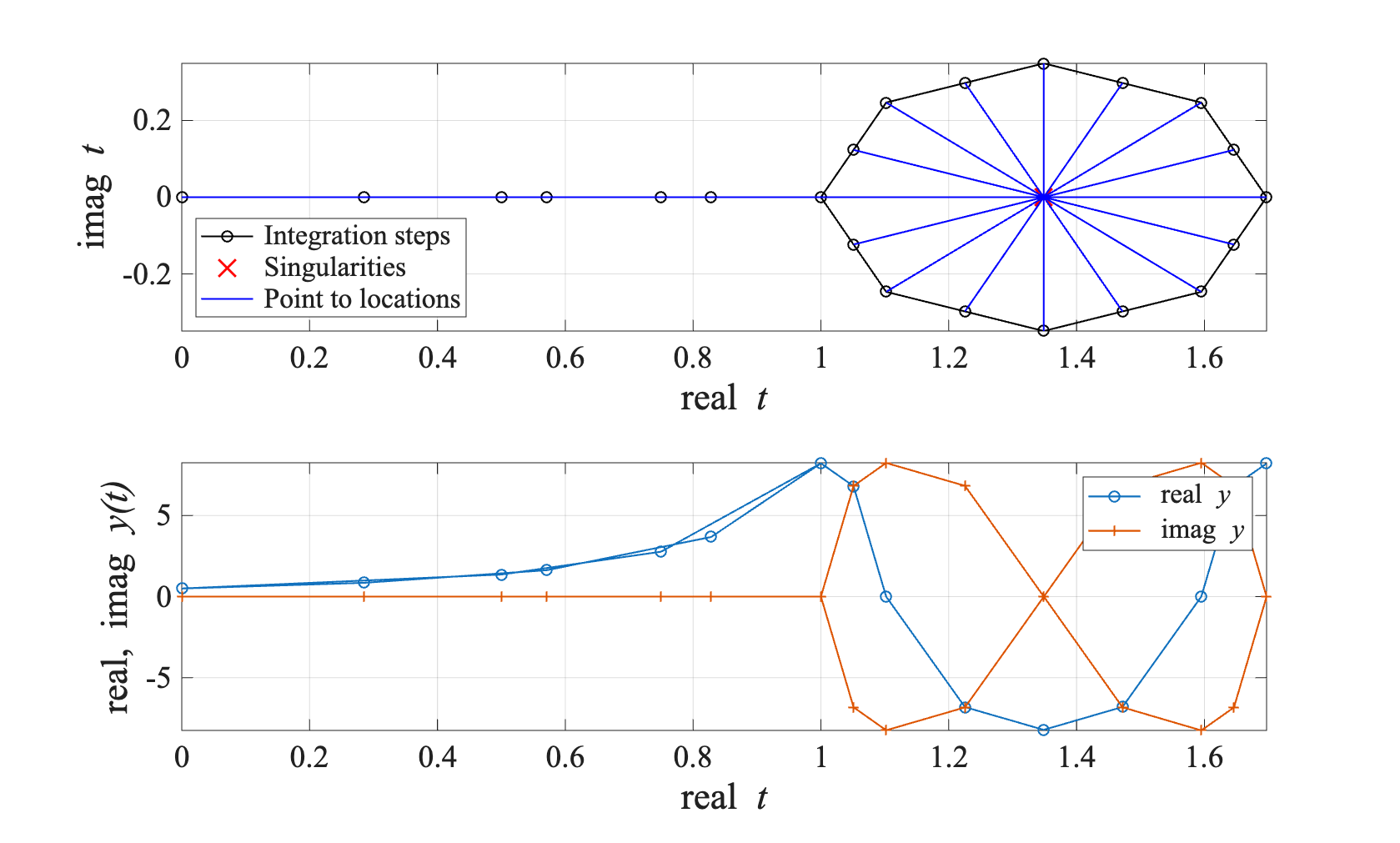}
  \end{center}
  \vspace{-20 pt}
  \caption{Circling a singularity remaining within the neighborhood of the pole.  The surfaces of the solution real and imaginary components are similar to those shown in Figure~\ref{fig:PainleveRealImag2}. As in Figure~\ref{fig:PainleveCirclePole}, the solution returns to its starting value.}
  \label{fig:PainleveTightPoleCircle}
\end{figure}

In the neighborhood of each singularity, the $\Real{y}$ and the $\Imag{y}$ surfaces are similar to those shown in Figure~\ref{fig:PainleveRealImag2}.  Instead of circling the pole on the real axis at $\Cmplx{1.3486}{0}$ following the offset circle shown in Figure~\ref{fig:PainleveCirclePole}, we follow a smaller circle centered at the pole as shown in Figure~\ref{fig:PainleveTightPoleCircle}.  The path shown in Figure~\ref{fig:PainleveCirclePole} intentionally  passes through the neighborhoods of poles at $3.1480 \pm 1.3306\, i$, leading to relatively large solution components in the vicinity of $\Real{t} \approx 2.8 - 3$ shown in the lower panel.  By keeping the path of integration in Figure~\ref{fig:PainleveTightPoleCircle} within the neighborhood of the pole at $\Cmplx{1.3486}{0}$, the solution components remain on surfaces very similar to those shown in Figure~\ref{fig:PainleveRealImag2}.

In Figure~\ref{fig:PainleveTightPoleCircle}, there are two unbounded peaks in the $\Real{y}$ component in the positive and negative $\Real{t}$ directions and two unbounded peaks in the $\Imag{y}$ component along the $45^{\rm o}$ lines.  There are two unbounded pits in the $\Real{y}$ component in the positive and negative $\Imag{t}$ directions and two unbounded pits in the $\Imag{y}$ components.  These are consistent with the surfaces shown in Figure~\ref{fig:PainleveRealImag2}.

An observer at the origin looking toward the singularity at $\Cmplx{1.3486}{0}$ in $\Real{y}$ faces a ridge rising into the clouds.  To the right and to the left of the singularity are unboundedly deep pits.  On the $45^{\rm o}$ lines to the right and to the left, the approach to the singularity is level with an increasing side-to-side gradient.  The same observer of the $\Imag{y}$ surface faces an unboundedly deep pit toward the right and an unboundedly high peak to the left.

Outside the neighborhood of a singularity, $y$ qualitatively resembles a linear combination of nearby singularities with weights approximated by the reciprocal of the distance from the observation point.  The influence of singularities more than twice the distance to the closest singularity is insignificant, as explained in Section~\ref{sec:RadiusOfConvergence}.

\begin{figure}[ht!]
  \begin{center}
  \includegraphics[width=0.9\linewidth]{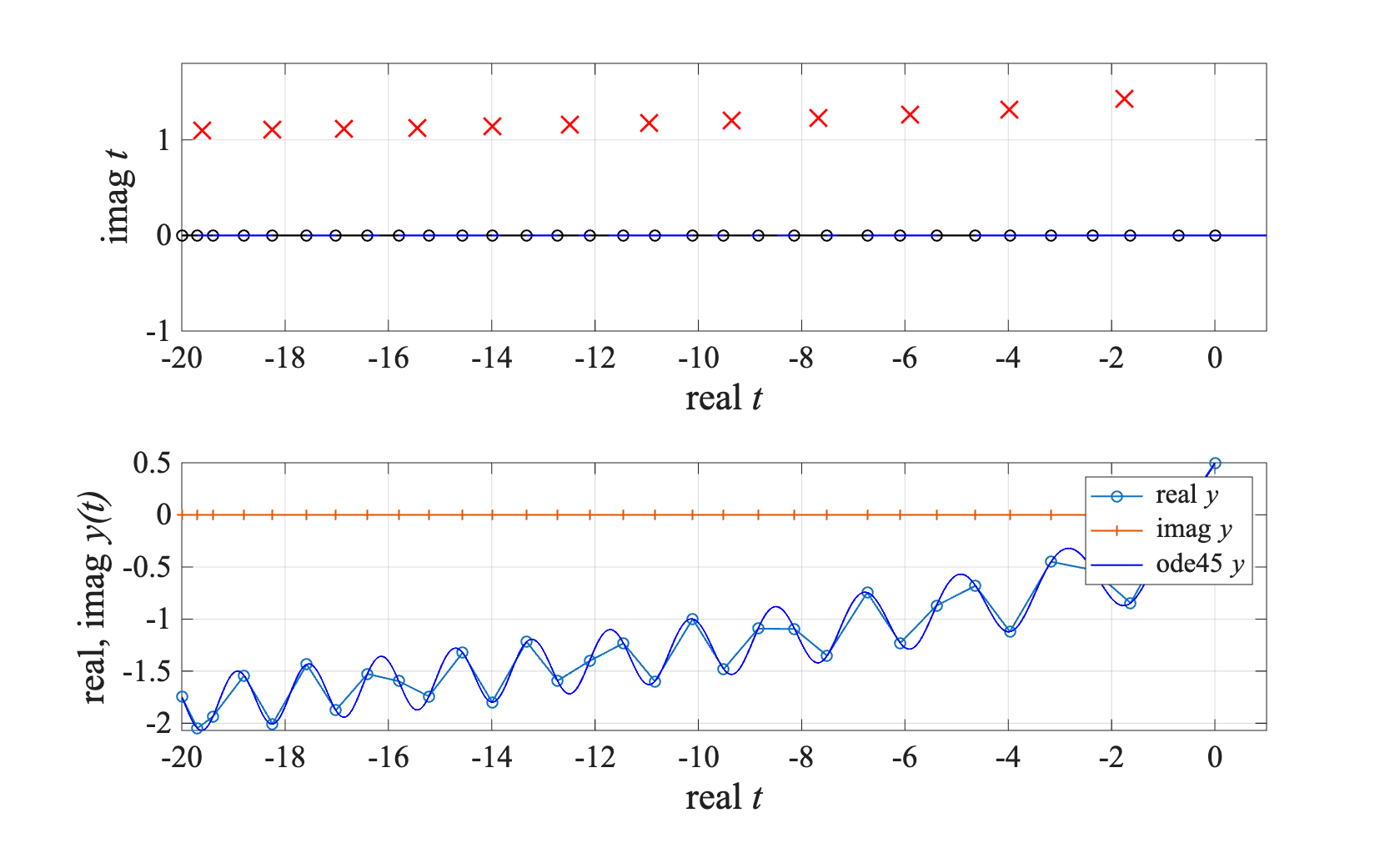}
  \end{center}
  \vspace{-20 pt}
  \caption{Integrating along the negative real $t$-axis.  The negative real $t$-axis is equidistant from conjugate pairs of singularities, so Algorithm 1 locates no singularities.  The singularities shown were located previously along a path of integration along $\Imag{t} \approx 0.75$. In the lower panel, $\Imag{y} = 0$ because the solution is real on the real axis.}
  \label{fig:PainleveNegativeRealAxis}
\end{figure}

A path of integration along the negative real $t$-axis encounters conjugate pairs of singularities with imaginary components slightly over one.  Figure~\ref{fig:PainleveNegativeRealAxis} shows only the singularities in the upper half plane.  In the lower panel, $\Imag{y} = 0$ because the true solution is real-valued on the real $t$ axis.  The computed real-valued solution tends downward, oscillating slightly under the influence of nearby singularities, like gently rolling hills.  Recall that the solution in the neighborhood of each singularity has unbounded peaks and pits, so the computed solution is relatively peaceful compared with the rough terrain closer to the singularities. As a landscape, it may have a very slight resemblance to some of the U.S. National Parks in Utah, with peaks on both sides of a long, narrow valley whose center is the negative $t$-axis.

\begin{figure}[ht!]
  \vspace{-10 pt}
  \begin{center}
  \includegraphics[width=0.9\linewidth]{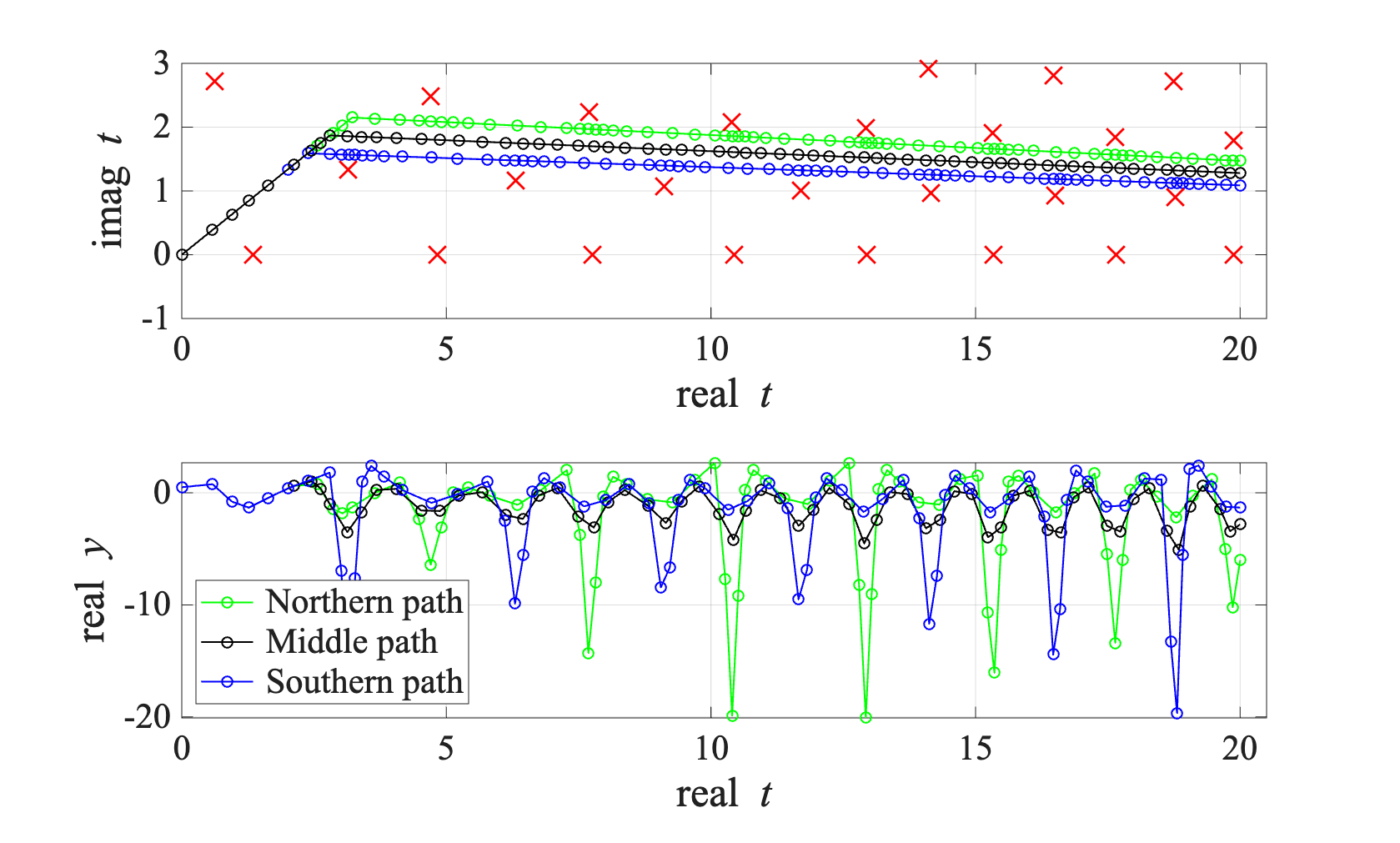}
  \end{center}
  \vspace{-20 pt}
  \caption{Upper panel shows three paths of integration between two rows of singularities in Region 1.  Lower panel shows $\Real{y}$ vs. $\Real{t}$ along the same three paths.}
  \label{fig:PainleveValley}
\end{figure}

Figure~\ref{fig:PainleveValley} explores the values of $\Real{y}$ in a valley between the second and third rows of singularities Region 1.  The path along the center of the ``valley'' (black) has variations less than six from a gently declining value, with oscillations from singularities alternating between above and below the path.  In contrast, paths closer to the singularities (green and blue) have larger oscillations.  The extreme values are negative because $\Real{y}$ becomes unboundedly negative south and north of each singularity.  The solution along the northern path appears more volatile than the solution along the southern path because the northern path passes closer to its nearby singularities than the southern path does.

\begin{figure}[ht!]
  \begin{center}
  \includegraphics[width=0.8\linewidth]{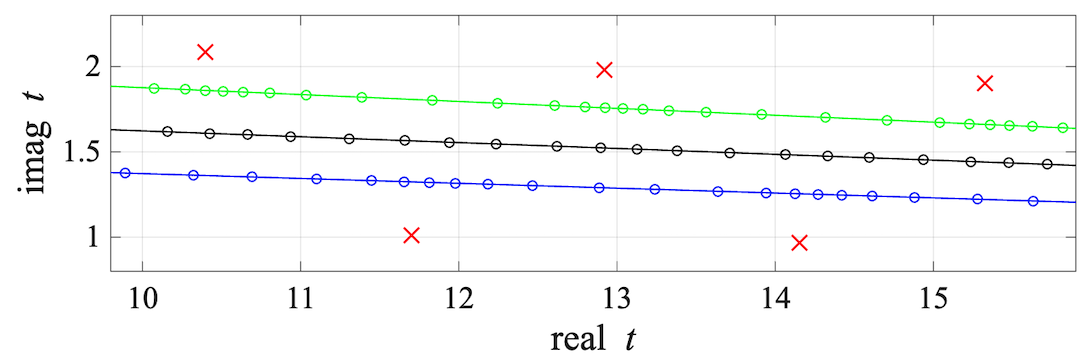}
  \end{center}
  \vspace{-15 pt}
  \caption{Detail in Figure~\ref{fig:PainleveValley} near $t = \Cmplx{13}{1.6}$ showing southern (blue), middle (black), and northern (green) paths. Step sizes are reduced near singularities because the {\tt odets} step size is nearly proportional to the series radius of convergence.}
  \label{fig:PainleveValley13}
\end{figure}

Figure~\ref{fig:PainleveValley13} zooms into Figure~\ref{fig:PainleveValley} near $t = \Cmplx{13}{1.6}$.  The northern (green) path and the southern (blue) path require smaller integration steps in the neighborhood of each singularity they pass because the {\tt odets} step size is nearly proportional to the series radius of convergence, determined by the distance to the closest singularity.  The middle (black) path can maintain accuracy with longer integration steps because the primary singularities on each side of the path of integration are further from the path. 

The solution shown in these graphs is not smooth because we plot the solution values only at the integration steps for economy.  With the 45-term Taylor series, the integration steps are quite long, and the details are lost. 
For smooth solutions, we can ask {\tt Odets} to evaluate the Taylor series for the solution at intermediate points.

\section{What is the value of the solution at $t = 20$?}
\label{sec:SolutionAtTwenty}

Given Problem~(\ref{Eq:IVP}), 
if we integrate in the direction of positive, real $t$, $y$ becomes unbounded, and we cannot advance beyond the pole of order two at $t \approx 1.3486$.  However, if we consider $t$ to be complex, the solution $y(t)$ is defined for all complex $t$, except for isolated singularities such as those shown in Figure~\ref{fig:PainleveLocations}.  Hence, $y(20)$ has a value.  Of course, $y(t)$ is real on the real $t$-axis, so $\Imag{y(20)} = 0$.

\begin{figure}[ht]
  \vspace{-10 pt}
  \begin{center}
  \includegraphics[width=0.9\linewidth]{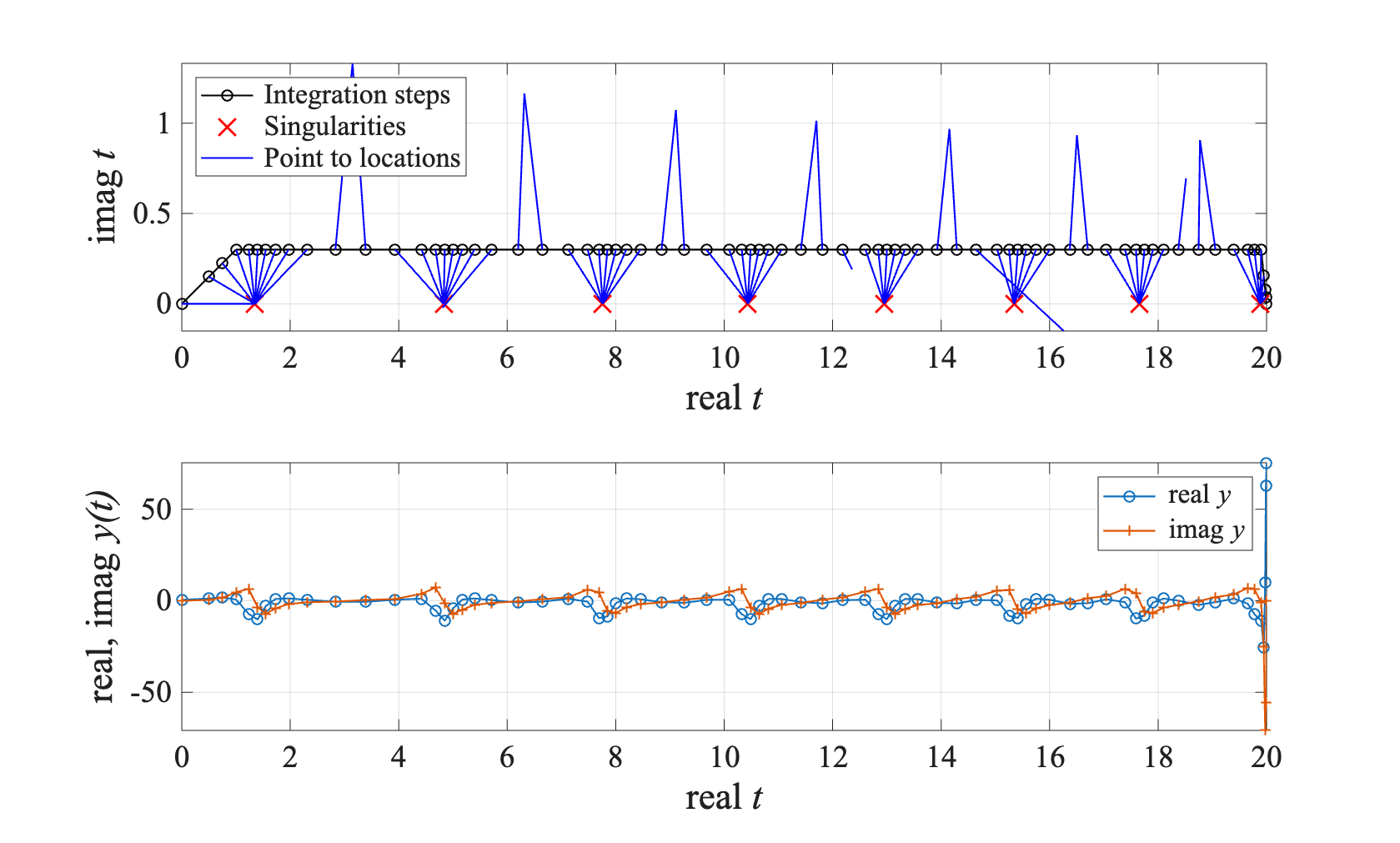}
  \end{center}
  \vspace{-20 pt}
  \caption{Solution $y(20) = 75.44599788$ by a ``long jump'' passing over eight singularities.}
  \label{fig:PainleveYAt20}
\end{figure}

\begin{figure}[ht]
  \vspace{-10 pt}
  \begin{center}
  \includegraphics[width=0.8\linewidth]{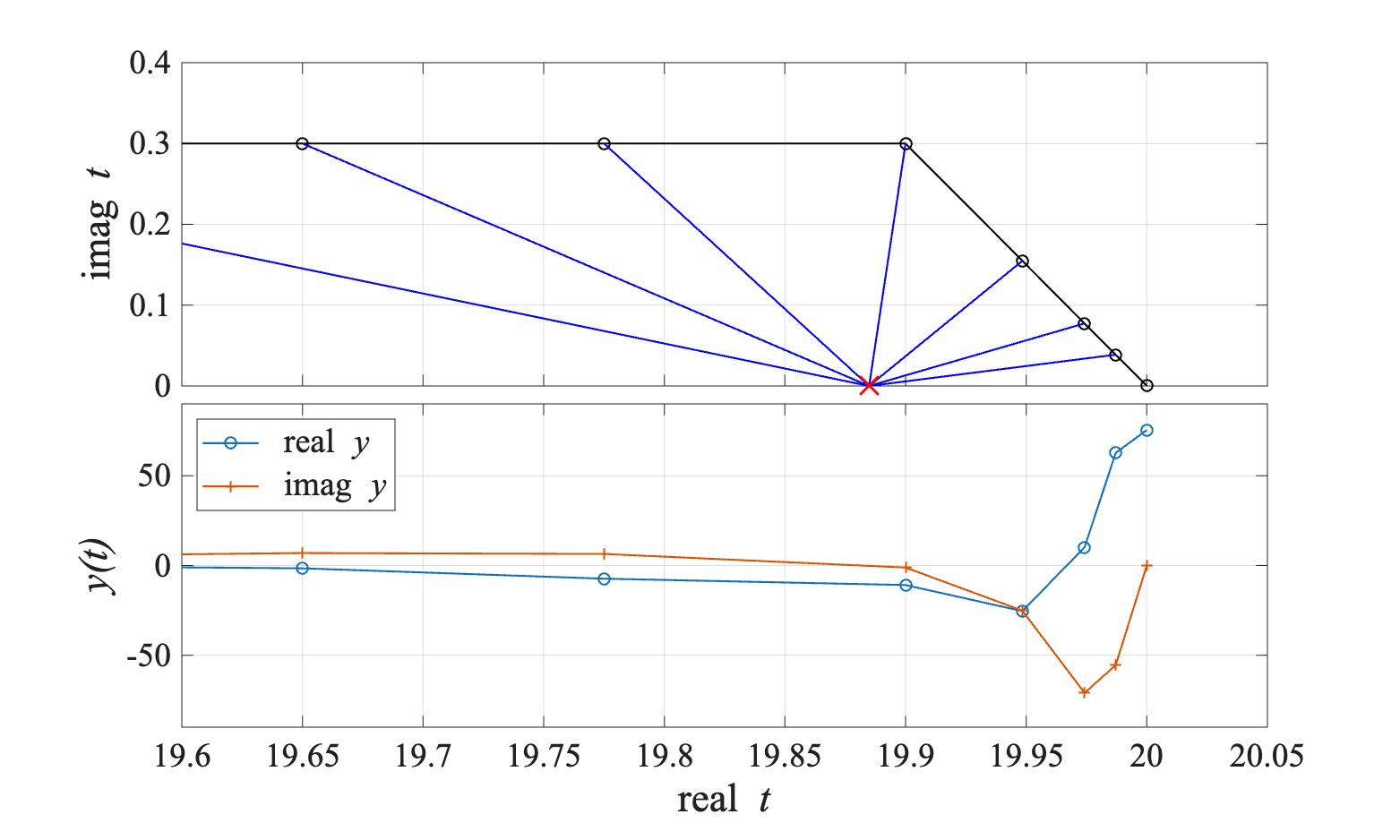}
  \end{center}
  \vspace{-20 pt}
  \caption{Solution $y(20) \approx 75.44599788$ by a ``long jump.'' Detail of Figure~\ref{fig:PainleveYAt20} near $t = 20$.}
  \label{fig:PainleveCircularLongJump}
\end{figure}

Let us compute $y(20)$.  We can follow the path of integration shown in Figure~\ref{fig:PainleveYAt20} to the points $t = \Cmplx{0.3}{0.3}$, $\Cmplx{20}{0.3}$, and then $\Cmplx{20}{0}$.  Figure~\ref{fig:PainleveYAt20} is similar to Figure~\ref{fig:PainleveRegion1A}, except that in Figure~\ref{fig:PainleveYAt20}, beyond $t = 20$, the path of integration returns to the $\Real{t}$-axis, as shown in detail in Figure~\ref{fig:PainleveCircularLongJump}.  We get 
  \[
    y(20) =  75.4459765085 + \num{1.4e-5}\, i,\ y'(20) = -1311.8172203538 - \num{3.6e-4}\, i\ .
  \]
The imaginary components of $y(20)$ and $y'(20)$ should be 0.0.  ${\rm{Real}}(y(20))$ and \\
$\Real{y'(20)}$ are larger than one might expect because there are unbounded peaks in $\Real{y}$ to the east and west of each singularity. We call the technique of integrating an ODE in the complex plane along a path to avoid real singularities a ``long jump.''

\begin{figure}[ht]
  \vspace{-10 pt}
  \begin{center}
  \includegraphics[width=\linewidth]{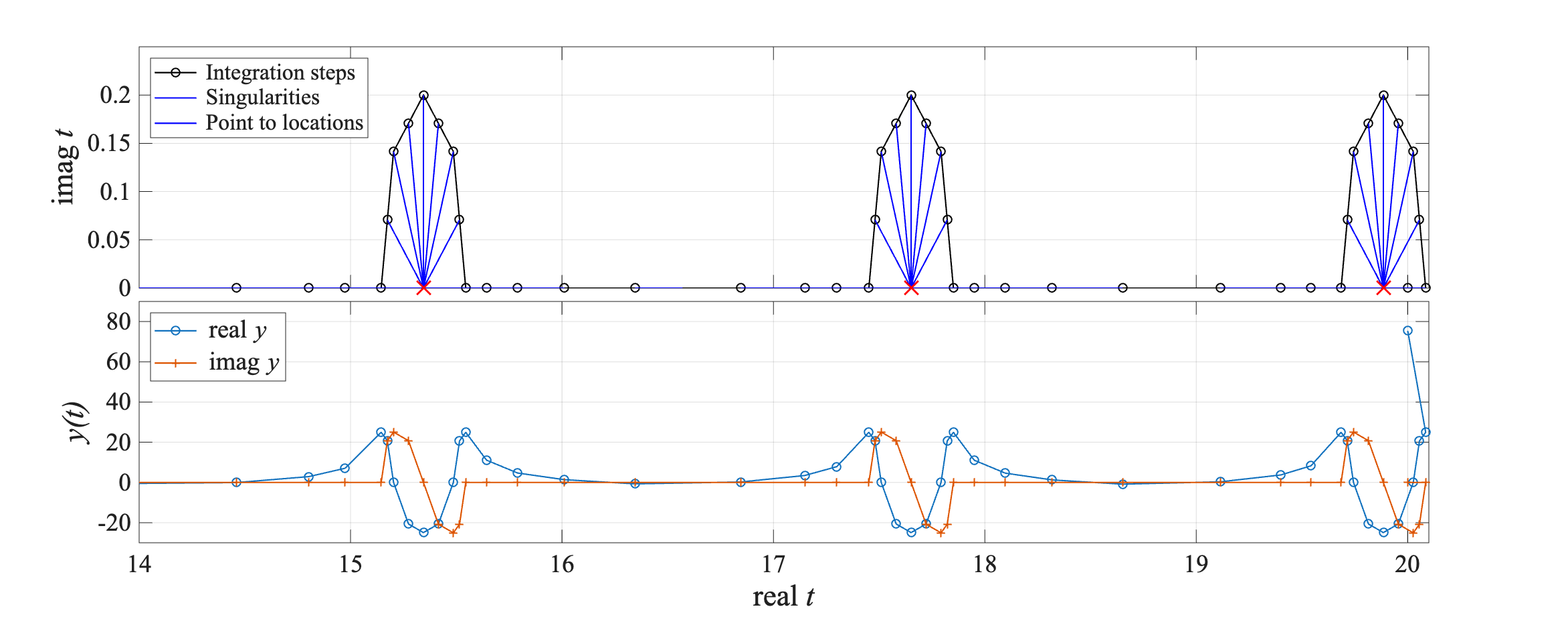}
  \end{center}
  \vspace{-15 pt}
  \caption{Solution $y(20) \approx  75.44599788$ by semi-circular ``pole vaulting,'' showing only the last three singularities.}
  \label{fig:PainleveCircularPoleVaults}
\end{figure}

Alternatively, if we want to know more about $y$ along the real $t$-axis, we can integrate along the real $t$-axis, except for taking small ``vaults'' along semicircles over each singularity as suggested in Figure~\ref{fig:PainleveCircularPoleVaults}.  In the upper panel of Figure~\ref{fig:PainleveCircularPoleVaults}, the vertical axis is four times the scale of the horizontal axis.  Using this technique, we get
  $$
    y(20) =  75.4459978841 + \num{1.2e-9}\, i, y'(20) = -1311.8177771832 - \num{3.3e-8}\, i\ .
  $$
The imaginary components have about twice as many accurate digits as for the ``long jump'' integration path.  Perhaps the real components are more accurate, too.

\begin{figure}[ht!]
  \vspace{-10 pt}
  \begin{center}
  \includegraphics[width=0.8\linewidth]{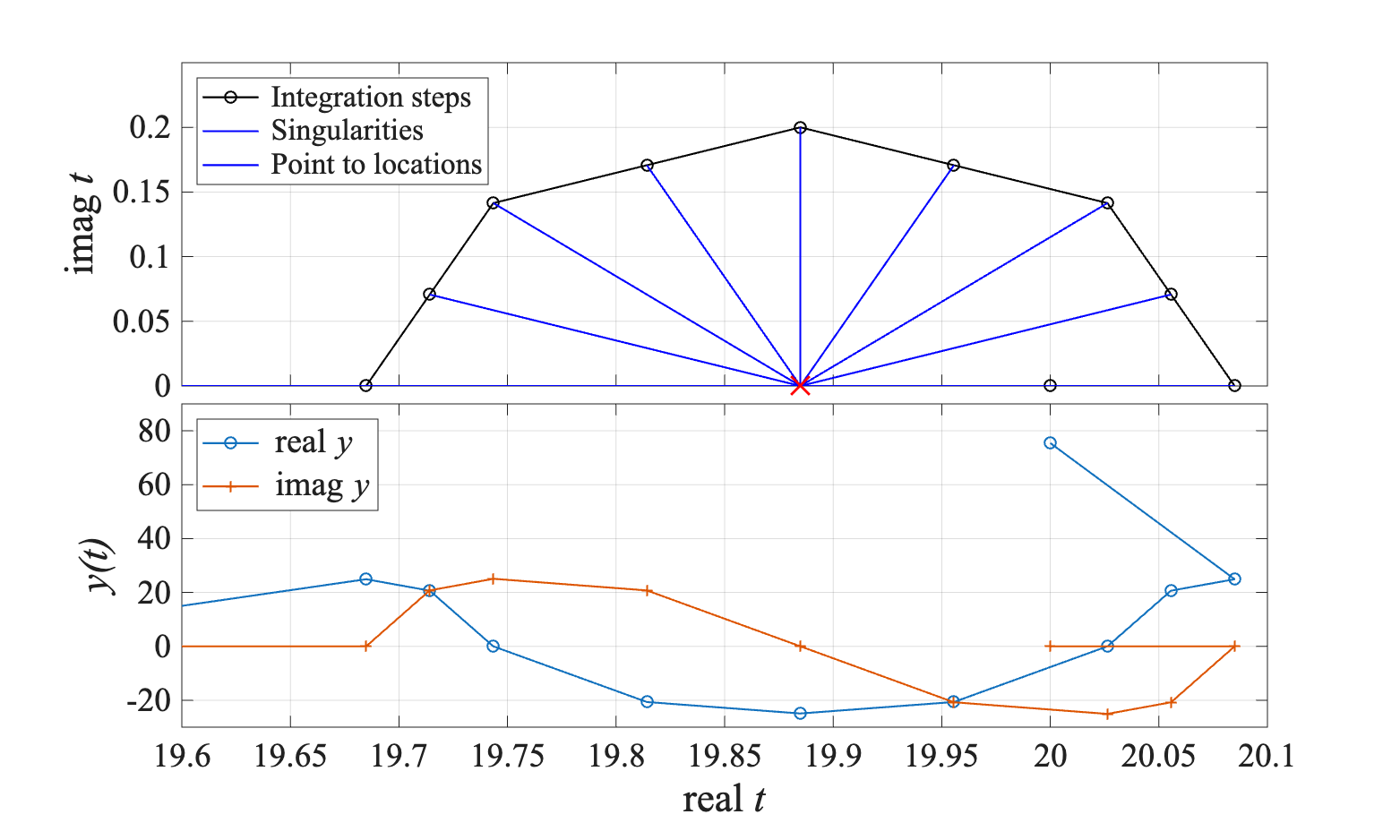}
  \end{center}
  \vspace{-20 pt}
  \caption{Solution $y(20) \approx  75.44599788$ by semi-circular ``pole vaulting.'' Detail near $t = 20$.}
  \label{fig:PainleveCircularPoleVaultsZoom}
\end{figure}

The lower panel of Figure~\ref{fig:PainleveCircularPoleVaults} offers a qualitative view of $y(t)$ along the real $t$-axis showing only the last three singularities.  Figure~\ref{fig:PainleveCircularPoleVaultsZoom} shows the detail near $t = 20.$  Between singularities, $\Real{y}$ is nearly zero, rising and falling sharply as $(a - t)^{-2}$ near each singularity.  We call the technique of integrating an ODE in the complex plane taking hops over each real singularity ``pole vaulting''~\cite{Corliss1980Integrating}.

\section{Conclusions}
\label{sec:Conclusions}

Figure~\ref{fig:PainleveLocations} illustrates the ability of \Alg:ref{alg:LocSing} to provide the numerical locations of hundreds of poles of the First \Pain Transcendent based on the work of Chang and Corliss~\cite{ChangCorliss1980ratio}. A similar technique in~\cite{ChangTaborWeissCorliss1981analytic} was applied to the Henon-Heiles System, whose singularities have two different fractional orders in a fractal structure forming a natural boundary.
We invite applications for which numerical maps of locations of singularities similar to those provided here can offer useful physical insights.

Sections~\ref{sec:ReproducableLocations} -~\ref{sec:SolutionAtTwenty} each pose and answer questions about the behavior of the surface traced by the First \Pain Transcendent employing singularity location information.  Taken together, these represent a sample of insights that may be disclosed by singularity location techniques.
  \begin{description}
    \item[\S 6] Are the locations we find reproducible, or do they depend on the integration path?  They are stable with respect to path.
    \item[\S 7] If we follow a closed integration path around a singularity, do we return to the same solution value?  Yes.
    \item[\S 8] Qualitatively, the complex-valued surface $y$ resembles a landscape of long, gentle-curving valleys with sharp peaks and pits on each side.
    \item[\S 9] Can we compute a meaningful value for the solution beyond several singularities on the real $t$-axis?  Yes, $y(20) \approx 75.44599788$.
  \end{description}

\section*{Acknowledgments}

Much of the work reported here was done in the early 1980s in collaboration with Y.F. Chang (now deceased) and Richard Kelley.
Thanks to John Pryce and Ned Nedialkov for many years of friendship and collaboration, for permission to use and experiment with {\tt odets}~\cite{NedialkovPryce2025}, and for comments that improved this work.  Special thanks are due to Karen Wurzel for many insightful suggestions and questions.

\section*{Funding}

The original work in the early 1980s was supported by the author's faculty salary from Marquette University.  The author returned to this project in retirement with no external support.

\section*{Conflict of interest}

The author has no conflict of interest associated with the work reported here.

\section*{AI tools}

The use of artificial intelligence tools in the preparation of this document was limited to spell checking, grammar checking, and sentence-level rewording for improved expression.

\bibliographystyle{elsarticle-num}
\bibliography{ODE}

@article{borisov2014paul,
	author = {Borisov, Alexey V and Kudryashov, Nikolay A},
	journal = {Regular and Chaotic Dynamics},
	pages = {1--19},
	publisher = {Springer},
	title = {Paul {P}ainlev{\'e} and His Contribution to Science},
	url = {https://link.springer.com/article/10.1134/S15603 5474010018},
	volume = {19},
	year = {2014}}

@article{Boutroux1913a,
	author = {Boutroux, M. P.},
	journal = {Annales Scientifiques de l'{\'{e}}.N.S. $3^e$ s{\'{e}}rie},
	pages = {255--375},
	publisher = {American Mathematical Society},
	title = {Recherches sur les transcendantes de {M}. {P}ainlev{\'{e}} et l'{\'{e}}tude asymptotique des {\'{e}}quations diff{\'{e}}rentielles du second ordre},
	url = {https://www.numdam.org/item/10.24033/asens.661.pdf},
	volume = {30},
	year = {1913},
    language = {French}}

@article{ChangTaborWeissCorliss1981analytic,
	author = {Chang, YF and Tabor, Michael and Weiss, J and Corliss, George F.},
	journal = {Physics Letters A},
	number = {4},
	pages = {211--213},
	publisher = {Elsevier},
	title = {On the Analytic Structure of the {H}enon-{H}eiles System},
	volume = {85},
	year = {1981}}

@article{ChangCorliss1980ratio,
  title={Ratio-like and Recurrence Relation Tests for Convergence of Series},
  author={Chang, YF and Corliss, George F.},
  journal={IMA Journal of Applied Mathematics},
  volume={25},
  number={4},
  pages={349--359},
  year={1980},
  publisher={Oxford University Press},
  doi={\ doi.org/10.1093/imamat/25.4.349}
}

@article{Corliss1980Integrating,
	author = {Corliss, George F.},
	journal = {Math. Comp.},
	number = {152},
	pages = {1181-1189},
	publisher = {American Mathematical Society},
	title = {Integrating {ODE}s in the Complex Plane - {P}ole Vaulting},
	url = {https://www.ams.org/journals/mcom/1980-35-152/S0025-5718 -1980-0583495-8/},
	volume = {35},
	year = {1980}}

@book{Fobes1963Calculus,
  title={Calculus and Analytic Geometry},
  author={Fobes, Melcher P. and Smyth, Ruth B.},
  number={v. 1},
  lccn={63007405},
  series={Calculus and Analytic Geometry},
  url={https://books.google.com/books?id=eiLYAAAAMAAJ},
  year={1963},
  publisher={Prentice-Hall}
}

@book{Hartman1964ordinary,
  title={Ordinary Differential Equations},
  author={Hartman, Philip},
  year={1964},
  city={New York},
  publisher={Wiley}
}

@book{Hille1976a,
	author = {Hille, Einar},
	city = {New York},
	publisher = {Wiley},
	title = {Ordinary Differential Equations in the Complex Domain},
	year = {1976}}

@book{Ince1927a,
	author = {Ince, Edward Lindsay},
	publisher = {Longmans, Green and Co. Ltd.},
	title = {Ordinary Differential Equations},
	year = {1927}}

@article{NedialkovPryce2025,
    author = {Nedialkov, Nedialko S. and Pryce, John D.},
    title = {Sub-{ODE}s Simplify {T}aylor Series Algorithms for Ordinary Differential Equations},
    journal = {SIAM Journal on Scientific Computing},
    volume = {47},
    number = {5},
    pages = {A2746-A2773},
    year = {2025},
    doi = {\ doi.org/10.1137/24M1716161}}

@inproceedings{Painleve1888lignes,
	author = {Painlev{\'e}, Paul},
	booktitle = {Annales de la Facult{\'e} des Sciences de Toulouse: Math{\'e}matiques},
	pages = {B1--B130},
	title = {Sur les lignes singuli{\`e}res des fonctions analytiques},
	volume = {2},
	year = {1888},
    URL = {https://archive.org/details/surleslignessin00paingoog},
    language = {French}}

@book{Painleve1895leccons,
	author = {Painlev{\'e}, Paul},
	publisher = {A. Hermann},
	title = {Le{\c{c}}ons sur l'Int{\'e}gration des {\'E}quations Diff{\'e}rentielles de la M{\'e}canique et Applications},
	url = {https://books.google.com/books?hl=en&lr=&id=-s1UAAAAYAAJ&oi=fnd&pg=PA3&ots=ZR8ZwuZzEf&sig=XkpKi5TLpT07werALf1qPo9tj70#v=onepage&q&f=false},
	year = {1895},
    language = {French}}

@article{Painleve1900memoire,
	author = {Painlev{\'e}, Paul},
	journal = {Bulletin de la Soci{\'e}t{\'e} math{\'e}matique de France},
	pages = {201--261},
	title = {M{\'e}moire sur les {\'e}quations diff{\'e}rentielles dont l'int{\'e}grale g{\'e}n{\'e}rale est uniforme},
	url = {https://www.numdam.org/item/10.24033/bsmf.633.pdf},
	volume = {28},
	year = {1900},
    language = {French}}

@article{Picard1880Propriete,
	author = {Picard, {\'E}mile},
	journal = {Bulletin des Sciences Math{\'e}matiques et Astronomiques},
	pages = {416-432},
	title = {Sur une Propriet{\'e} des Fonctions Uniformes d'une Variable Li{\'e}es par une Relation Alg{\'e}brique et sur une Classe d'{\'e}quations Diff{\'e}rentiels},
	volume = {4 ($2^{{\rm e}}$ s{\'e}rie)},
	year = {1880},
    url = {https://www.numdam.org/item/BSMA_1880_2_4_1_416_0/},
    language = {French}}

@article{Poincare1885Comptes,
	author = {Poincar{\'e}, Henri},
	journal = {Acta Mathematica},
	pages = {1--32},
	title = {{Sur un th\'eor\`eme de M. Fuchs}},
	volume = {7},
    number = {1},
	year = {1885},
    url = {https://henripoincarepapers.univ-lorraine.fr/chp/hp-pdf/hp1885am.pdf},
    language = {French}}

\end{document}